\documentclass[a4paper,10pt]{article}
\usepackage{makeidx}
\usepackage[french,english]{babel}
\usepackage[latin1]{inputenc}
\usepackage{amsfonts}
\usepackage{bm}
\usepackage{amssymb}
\usepackage{latexsym}
\usepackage{amsmath}
\usepackage{amscd}
\usepackage{amsthm}
\usepackage{hyperref}
\usepackage{rotating}
\usepackage{graphicx}
\usepackage{pifont}
\usepackage{curves}
\usepackage{fancyhdr}
\usepackage{epsfig}
\usepackage{pstricks}
\usepackage{pst-tree}
\usepackage{subfigure}
\usepackage{epic}
\usepackage{alltt}
\usepackage[all]{xy}
\usepackage{appendix}

\newcommand{\Co}{\mathcal{C}}
\newcommand{\D}{\mathbb{D}}

\newcommand{\Card}{\mbox{Card}}

\newcommand{\N}{\mathbb{N}}

\newcommand{\R}{\mathbb{R}}

\renewcommand{\P}{\mathbb{P}}
\newcommand{\E}{\mathbb{E}}
\newcommand{\ind}{{\mathchoice {\rm 1\mskip-4mu l} {\rm 1\mskip-4mu l}
{\rm 1\mskip-4.5mu l} {\rm 1\mskip-5mu l}}}

\theoremstyle{plain}
\newtheorem{theorem}{Theorem}[section]
\newtheorem{proposition}[theorem]{Proposition}
\newtheorem{lemma}[theorem]{Lemma}
\newtheorem{corollary}[theorem]{Corollary}

\theoremstyle{definition}
\newtheorem{definition}[theorem]{Definition}
\newtheorem{hyp}[theorem]{Assumption}

\title{Daphnias: from the individual based model to the large population equation}

\author{J.A.J.(Hans) Metz \and Viet Chi Tran}

\author{J.A.J.(Hans) Metz and Viet Chi Tran\footnote{J.A.J. Metz: Mathematical Institute \& Institute of Biology \& NCB Naturalis, Leiden; Niels Bohrweg 1, 2333 CA, Leiden, Netherlands \& EEP, IIASA, Laxenburg, Austria. V.C. Tran: Laboratoire Paul Painlevé, UFR de Mathématiques, UMR CNRS 8524; Cité Scientifique, 59 655 Villeneuve d'Ascq Cédex, France. \& CMAP, Ecole Polytechnique.}}

\begin{document}

\maketitle

\begin{center} {\small {\it Dedicated to Odo Diekmann on the occasion of his 65th birthday}}
\end{center}
\bigskip

\begin{abstract}The class of deterministic 'Daphnia' models treated by Diekmann et~al. (J Math Biol   61: 277--318, 2010) has a long history going back to Nisbet and Gurney (Theor Pop Biol 23: 114--135, 1983)  and Diekmann et~al. (Nieuw Archief voor Wiskunde  4: 82--109, 1984).  In this note, we formulate the individual based models (IBM) supposedly underlying those  deterministic models. The models treat the interaction between a general size-structured consumer population ('Daphnia') and an unstructured resource ('algae'). The discrete, size and age-structured Daphnia population changes through births and deaths of its individuals and throught their aging and growth. The birth and death rates depend on the sizes of the individuals and on the concentration of the algae. The latter is supposed to be a continuous variable with a deterministic dynamics that depends on the Daphnia population. In this model setting we prove that when the Daphnia population is large, the stochastic differential equation describing the IBM can be approximated by the delay equation featured in  (Diekmann et~al., l.c.).
\end{abstract}


\section{Introduction}

The theory of physiologically structured populations as developed in works by Diekmann and co-authors \cite{MetzDiekmann,Diekmannetal1998,Diekmannetal2001,Diekmannetal2003,diekmanngyllenbergmetznakaokaderoos} derives its motivation from its supposed ability to link population level phenomena to specific mechanisms in and around individual organisms. Yet, those individuals do not figure as such in the models, which treat the spatial concentrations of those individuals as continua. The supposition has always been that the proposed deterministic frameworks would in principle be derivable from individual-based stochastic models  (c.f. \cite{MetzRoos}), but full derivations of this sort so far have only been done for the special cases of finite i-state models  (i from individual; leading to ODEs in the large number limit, e.g.\cite{Kurtz1970,Kurtz1981}) and age-based models (e.g. \cite{oelschlager,chithese,trangdesdev,ferrieretran,jagersklebaner,jagersklebaner2}). Moreover, it is possible to reinterpret the results in \cite{champagnatferrieremeleard} as pertaining to general structured models with only jump transitions.
Although some age-structured models are phrased in terms of i-state variables that change with time in a fixed manner, uninfluenced by the environment, these can truly represent but a small fraction of the rich variety of physiological mechanisms seen in nature (e.g. \cite[chapters I and III]{DiekmannMetz2010,MetzDiekmann}). For example, for the majority of species, size with a growth rate that depends on the environment is a far more important determinant of an individual's population dynamical behaviour than is age. Except in the physiologically well buffered homeotherms (on which we anthropocentrically are inclined to focus) reproduction tends to be under a strong influence of past food availabilities, with reproductive size being reached far earlier when an individual has encountered good than when it has encountered bad feeding conditions.  In, for example, arthropods or fishes, maturing to the reproductive states often largely depends on how much food they have encountered but little on their age. However, the toolbox for proving the appropriate law of large number results that was developed in \cite{fourniermeleard,chithese,trangdesdev,ferrieretran} needs to be extended a bit to deal with size-based models, referred to as 'Daphnia' models by  \cite{diekmanngyllenbergmetznakaokaderoos}, let alone to the even more general models considered in \cite{MetzDiekmann,MetzRoos,Diekmannetal1998,Diekmannetal2001,Diekmannetal2003}. In this paper we, as a birthday present to Odo Diekmann, intend to remedy this lack,  although so far only for the 'Daphnia' models that appear as central example in his work  (c.f. \cite{diekmanngyllenbergmetznakaokaderoos}).

Diekmann et~al. \cite{diekmanngyllenbergmetznakaokaderoos} lay out a general framework for investigating the interaction between a general age or size-structured consumer population (referred to as 'Daphnia') and an unstructured resource (referred to as 'algae'), a class of problems special forms of which were first considered in \cite{NisbetGurney,gurneynisbet} and  \cite{DiekmannMetzKooijmanHeymans1984,deroosmetzeversleipoldt}. Diekmann et~al. (l.c.) show that stability properties and bifurcation phenomena can be understood in terms of solutions of a system of two delay equations that are analysed using results in \cite{diekmanngettogyllenberg,DiekmannGyllenberg2012}. In this note, we derive their model from a microscopic description: starting from a stochastic, age- and size-structured, individual based model (IBM) for the Daphnia population we recover the equations of \cite{diekmanngyllenbergmetznakaokaderoos}.\\
In our study, the Daphnia population is discrete and stochastic while the algal population is continuous.  We treat the Daphnia population as a point measure on a space spanned by size and age and from this platform generalize the microscopic construction given in \cite{trangdesdev} to arrive at a stochastic differential equation (SDE in the sequel) driven by Poisson point processes that gives a pathwise description of the population evolution. The main novelty is that the individual growth rate also depends on the population as a whole, instead of only on an individuals' characteristics, since the algal concentration and thus each individuals' resource access is affected by the entire Daphnia population.\\
We then, in the wake of \cite{oelschlager,chithese,trangdesdev,ferrieretran,jagersklebaner,jagersklebaner2},  provide a law of large numbers that allows approximating the dynamics of the individual-based size- and age-structured model by means of a partial differential equation (PDE) when the volume containing the Daphnia and algae as well as the population sizes are large. In \cite{meleardtran} similar limit theorems are used as a basis for obtaining adaptive dynamics approximations for the evolution of hereditary traits in age-structured populations, while \cite{durinxmetzmeszena} gives informed conjectures about such approximations for the more encompassing model classes treated in \cite{Diekmannetal2003}.\\
Finally we study the  limiting equations to recover the macroscopic (deterministic) system considered in \cite{diekmanngyllenbergmetznakaokaderoos}.\\

\noindent \textbf{Notation:} We will denote the set of finite measures on $\R_+^2$ endowed with the weak convergence topology as $\mathcal{M}_F(\R_+^2)$. For $\mu\in \mathcal{M}_F(\R_+^2)$ and a real measurable bounded function $f$ on this space, we set $\langle \mu, f\rangle=\int_{\R_+^2} f(\xi,a)\mu(d\xi,da)$. The set of bounded  real functions of class $\Co^1$ on $E$  with bounded derivatives is denoted as $\Co^1_{\rm{b}}(E,\R)$.\\
If $X$ is a process indexed by time, then we will denote  the value of $X$ at time $t$ as either  $X(t)$ or $X_t$ (the latter to avoid formulas becoming cluttered with too many brackets).

\section{Individual-based Daphnia model}

\subsection{Model specification}
Our Daphnia population consists of discrete individuals living in continuous time, differentiated by an age $a\in \R_+$ and a size $\xi\in \R_+$. Individuals are given labels $i\in \N^*=\{1,2,\dots\}$, with the individuals present in the population at $t=0$ bearing numbers in an order opposite to that of their ages,  followed by the individuals born after $t=0$  in the order of their appearance in the population. We denote as $I_t\in \N^*$ the total number of individuals that were present at time $0$ or were born between time 0 and time $t$. $V_t\subset \{1,\dots,I_t\}$ denotes the set of individual alive at time $t$. We can then represent the Daphnia population as a point measure on the state space $\N^*\times \R_+\times \R_+$:


\medskip
\begin{equation}
\bar{Z}_t=\sum_{i\in V_t}\delta_{(i,\xi^i_t,a^i_t)},
\end{equation}
where $\xi^i_t$ and $a^i_t$ are the size and age of individual $i$ at time $t$. We denote as $Z_t(d\xi,da)=\bar{Z}_t(\N^* \times d\xi\times da)$ the marginal measure of $\bar{Z}_t$ on $\R_+^2$; $V_t$ equals the support of the marginal measure of $\bar{Z}_t$ on $\N^*$. Size and age are related as follows.\\
An individual's age is equal to $t-t_0$ where $t$ is the current time and $t_0$ its birth time. Individuals grow up in an environment that is characterized at time $t$ by the algal concentration $S(t)$. All individuals are born with the same size $\xi_0$ (to keep things simple) and an individual with size $\xi$ at time $t$ grows at speed $g(\xi,S(t))$, so that the size of an individual aged $a$ born at time $t_0$ is:
\begin{equation}
\xi(a; t_0)=\xi_0+\int_0^a g(\xi(\alpha; t_0) ,S(t_0+\alpha))d\alpha.\label{growth}
\end{equation}
The growth rates depend on the other individuals in the population. The latter point is the novelty of the IBM presented here. (We could have taken hereditary traits as well as further age-like i-state variables on board as in \cite{trangdesdev}, but have decided not do this in order not to unduly complicate the story).\\
Reproduction is asexual. The birth and death rates of an individual with size $\xi$ and age $a$ at time $t$ are $\beta(\xi,a,S(t))$ respectively $\mu(\xi,a,S(t))$. \\
An individual with size $\xi$ depletes the food density at rate ${1\over K}\gamma(\xi,S(t))$ and the food density replenishes with rate $f(S(t))$ so that the food concentration $S(t)$ evolves according to
\begin{align}
\frac{dS}{dt}(t)= & f(S_t)-{1\over K}\int\gamma(\xi,S_t)Z_t(d\xi,da)=  f(S_t)-{1\over K}\sum_{i\in V_t}\gamma(\xi^i_t,S_t).\label{food}
\end{align}
For the biological justification think of the Daphnia population as living in a container of size $K$, so that the Daphnia density is ${|Z_t|\over K}$, with $|Z_t|:=\Card(V_t)=\big\langle Z_t,1\big\rangle$ the number of individuals.

\noindent
\begin{hyp}\label{hypotheses-taux}
In the sequel, we assume that the growth speed $g(\xi,S)$, the rates $f(S)$ and $\gamma(\xi,S)$ are continuous bounded functions and that\\
(i) the birth rate $\beta(\xi,a,S)$ is piecewise continuous and bounded by $\bar{\beta}$. \\
(ii) the death rate $\mu(\xi,a,S)$ is continuous and there exists a function $\underline{\mu}(a)$ and a constant $A\in (0,+\infty]$ such that $\forall (\xi,a,S)\in \R_+^3,$ $\mu(\xi,a,S)\geq \underline{\mu}(a)$ and $\int_0^A \underline{\mu}(a)da=+\infty$.\\
(iii) $g$ is Lipschitz continuous with respect to $\xi$, uniformly in $S$ on compact intervals of $\R_+$, and bounded by $\bar{g}$.\\
(iv) $f$ and $\gamma$ are uniformly Lipschitz continuous with respect to $S$ uniformly in $\xi$ on compact intervals of $\R_+$.
\end{hyp}

\noindent The assumption (i) on the birth rate ensures that in a short time interval a single individual can not beget too many young: intervals between births are stochastically lower bounded by exponential random variables with rate $\bar{\beta}$. The assumption (ii) on the death rate implies that individuals a.s. have lifetimes bounded by $A$. Finally the assumptions
 (iii) and (iv) ensure that there exist unique continuous solutions to \eqref{growth} and \eqref{food} as long as the number of individuals $|Z_t|$ remains finite, the latter being guaranteed by the fact that $|Z_t|$ is stochastically bounded by a pure birth process with birth rate $\bar{\beta}$.


For $\Phi\in \Co^1_{\rm{b}}(\R^2,\R)$ and $\varphi\in \Co^1_{\rm{b}}(\R_+^2,\R)$, we denote by $\Phi_{\phi}$ the function on $\mathcal{M}_F(\R_+^2)\times \R_+$ defined by $\Phi_{\varphi}(Z,S)=\Phi(\langle Z,\varphi\rangle,S)$. From the description of the population dynamics, it follows that the process $(Z(t),S(t))_{t\in \R_+}$ is characterized by the infinitesimal generator $L$ operating on the functions $\Phi_{\varphi}$:
\begin{align}
L\Phi_{\varphi}(Z,S)=\ & \partial_1 \Phi_{\varphi}(Z,S)\ \big\langle Z,g(.,S),\partial_\xi f(.,.)+\partial_a f(.,.) \big\rangle\nonumber\\
 + &\Big\langle Z, \mu(.,.,S)\Big(\Phi\big(\langle Z,\varphi\rangle-\varphi(.,.),S\big)-\Phi_{\varphi}\big(Z,S\big)\Big) \Big\rangle\nonumber\\
 +&  \Big\langle Z, \beta(.,.,S)\Big(\Phi\big(\langle Z,\varphi\rangle+\varphi(\xi_0,0),S\big)-\Phi_{\varphi}\big(Z,S\big)\Big)\Big\rangle\nonumber\\
  +& \partial_2 \Phi_{\varphi}(Z,S)\Big(f(S)-\big\langle Z,  {1 \over K}\gamma(.,S)\big\rangle\Big)\label{generateurinfinitesimal}
\end{align}The first term describes the aging and growth of the living individuals of the population. The second and third terms represent the demography of the population (deaths and births). The fourth term corresponds to the variation of the food. \\

In the next Subsection we introduce the pathwise construction of an IBM with the described dynamics and give an SDE driven by a Poisson point process that admits \eqref{generateurinfinitesimal} as infinitesimal generator. This is useful for simulations and for deriving moment conditions and large population approximations (e.g. \cite{champagnatferrieremeleard}). The evolution is piecewise deterministic: The size of the population is modified at birth or death events. Between these, conditionally on the structure of the population after the last event, the growth of the individuals and the food dynamics are deterministic.

\subsection{Construction of the IBM and a useful SDE}

Let us start with some heuristics. Consider at time $t$ a population given by $\bar{Z}_t$ and food concentration $S(t)$. If no birth or death event occurs between time $t$ and $t+s$, then $V_{t+s}=V_t$. The sizes $\xi^i(t+s)$ for $i\in V_t$ at time $t+s$ and the food concentration $S(t+s)$ are obtained by solving:
\begin{align}
& \xi^i(t+s)=\xi^i(t)+\int_t^{t+s} g(\xi^i(\tau),S(\tau))d\tau\label{eq1} \\
& S(t+s)=S(t)+\int_t^{t+s} \Big(f(S(\tau))-{1 \over K}\sum_{i\in V_t}\gamma(\xi^i(\tau),S(\tau))\Big)d\tau.\nonumber
\end{align}
Under Assumptions \ref{hypotheses-taux}, this system has a unique solution, which we denote as
$(\Xi^i(t+s ; t,\bar{Z}_t,S_t),\Sigma(t+s ; t,\bar{Z}_t,S_t) \ ;\  s\in \R_+,i\in V_t)$.
We will denote the coordinate of the  flow corresponding to \eqref{eq1} for an individual with initial condition $\xi$ as $\Xi(t+. ; t,\xi,\bar{Z}_t,S_t)$, so that in particular, $\Xi^i(t+s ; t,\bar{Z}_t,S_t)=\Xi(t+s ; t,\xi^i_t,\bar{Z}_t,S_t)$. In the sequel, we will also use that, if no births or deaths occur, for all $0\leq s\leq t$:
\begin{equation}
\Xi(t ; s,\Xi(s ; 0,\xi_0,\bar{Z}_0,S_0),\bar{Z}_t,S_t)=\Xi(t ; 0,\xi_0,\bar{Z}_0,S_0).\label{composition-flots}
\end{equation}
After a birth or a death the process is restarted with appropriately adapted initial conditions at that instant.

\begin{figure}[ht]
  \begin{center}
    \begin{picture}(250,180)(-20,-10)
      \put(0,0){\vector(1,0){250}}
      \put(0,0){\vector(0,1){160}}
      \put(-16, 160){size}
      \put(245, -10){time}
      \dashline{4}(130,0)(130,165)
      \put(129, -10){$s_2$}
      \dashline{4}(50,0)(50,165)
      \put(49, -10){$s_1$}
      \dashline{4}(240,0)(240,165)
      \put(239, -10){$t$}
      \put(-10, 38){$\xi'$}
      \put(242, 143){$\xi_1$}
      \put(242, 114){$\xi_2$}
      \put(242, 130){$\xi_3$}
      \put(240,135){\circle*{5}}
      \put(126,48){$\times$}
      \put(0,-10){0}
      \put(-10,0){$\xi_0$}
      \linethickness{1pt}
      \qbezier(0,40)(20,60)(50,75)
      \qbezier[40](50,75)(140,130)(240,140)
      \qbezier(50,75)(80,90)(130,100)
      \qbezier[50](130,100)(180,110)(240,115)
      \qbezier(130,100)(180,130)(240,135)
      \qbezier(50,0)(83,40)(129,50)
    \end{picture}
  \end{center}
  \vspace{-0.3cm}
  \caption{{\small\textit{At time 0 there is a single particle of size $\xi'$ and age $a'$, which is expected to have size $\xi_1=\Xi(t;0,\xi',\bar{Z}_0,S_0)$ at time $t$, where $\bar{Z}_0=\delta_{(1,\xi',a')}$. At time $s_1$, a second particle is born. Just before birth, the population is $\bar{Z}_{s_{1-}}=\delta_{(1,\Xi(s_1;0,\xi',\bar{Z}_0,S_0),a'+s_1)}$. After the birth, the size expected for the first particle at time $t$ is changed from $\xi_1$ to $\xi_2=\Xi(t; s_1,\Xi(s_1;0,\xi',\bar{Z}_0,S_0),\bar{Z}_{s_{1-}}+\delta_{(2,\xi_0,0)},S_{s_1-})$ since there will be less resources for the two particles. At time $s_2$, the second particle dies. The size expected for the first particle at $t$ is changed again from $\xi_2$ to $\xi_3$ as there is now more resources for the first particles's growth.}}} \label{Fig1}
\end{figure}
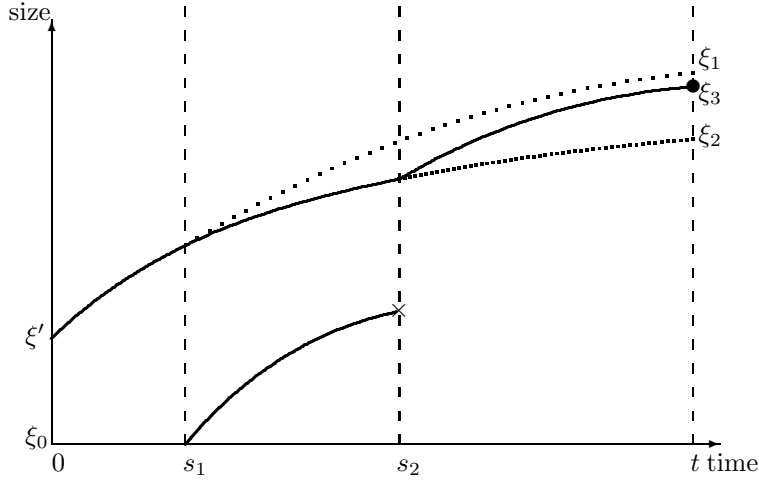

The above description suggests a simple direct way for simulating the IBM. Starting from a birth or death event, first generate a standard exponentially distributed random number $\tau$, and then run the differential equations for the states of all Daphnia individuals and for the algae. Simultaneously integrate the sum of the birth and death rates of the Daphnia, starting from zero.  When this integral reaches $\tau$, one of the Daphnia dies or gives birth.
Which individual is the culprit and whether the event is a birth or a death is then decided from a single multinomial draw with probabilities proportional to the contributions of all the different events to the total event rate at that time.\\
For the SDE that describes the process $(\bar{Z}(t),S(t))_{t\in \R_+}$, we proceed as in \cite{trangdesdev}, following a construction introduced by \cite{fourniermeleard} for the case without age or size, while accounting for the additional difficulty that the growth rate now depends on the rest of the population. To this end we again use that between two birth or death events the evolution of the population, conditionally on its state at the last event, is deterministic.
As the integral form in which we present the SDE looks back at the end result of the events happening over a time interval, we also have to look in a retrospective manner  at the resetting of the initial conditions at the moments that a birth or death occurs.
Assume that the initial condition $\bar{Z}_0$ and the initial food concentration $S_0$ are given. The idea is that to construct the population at time $t$, we can proceed as follows:
\begin{itemize}
\item If no event happens during $[0,t]$, then it is sufficient to consider the predicted sizes $(\Xi^i(t ; 0,\bar{Z}_0,S_0)\ ; \ i\in V_0)$ of individuals at $t$. The population at time $t$ is
$$\bar{Z}_t=\sum_{i\in V_0}\delta_{(i,\Xi^i(t ; 0,\bar{Z}_0,S_0),a^i_0+t)}$$and the food concentration is $S_t=\Sigma(t ; 0,\bar{Z}_0,S_0)$.
\item If a birth event occurs at time $s\in [0,t]$, then $V_s=V_{s_-}\cup \{I_{s_-}+1\}$, where $I_{s_-}$ is the number of labels already used so that the new individual gets labelled with the first available number $j=I_{s_-}+1$. The predicted sizes at time $t$, $(\Xi^i(t ; s,\bar{Z}_{s_-},S_{s_-})\ ; \ i\in V_{s_-})$, that we had for the individuals $i\in V_{s_-}$ before the event on the supposition that it were not to occur, are replaced by $(\Xi^i(t ; s,\bar{Z}_{s_-}+\delta_{(j,\xi_0,0)},S_{s_-})\ ; \ i\in V_{s_-})$; see Fig. \ref{Fig1}. Moreover for the new individual with label $j$, we add a Dirac mass at $(j,\Xi^{j}(t ; s, \bar{Z}_{s_-}+\delta_{(j,\xi_0,0)},S_{s_-}),t-s)$.
\item If the individual $j\in V_{s_-}$ dies at time $s\in [0,t]$, then $V_s=V_{s_-}\setminus \{j\}$ and the predicted sizes at time $t$, $(\Xi^i(t ; s,\bar{Z}_{s_-},S_{s_-})\ ; \ i\in V_{s})$, that we had for the individuals $i\in V_{s}$ before the event on the supposition that it were not to occur, are replaced by $(\Xi^i(t ; s,\bar{Z}_{s_-}-\delta_{(j,\xi^j_{s_-},a^j_{s_-})},S_{s_-})\ ; \ i\in V_{s})$; see Fig. \ref{Fig1}. Additionally, we delete the Dirac mass at $(j,\Xi^j(t ; s,\bar{Z}_{s_-},S_{s_-}),a^j_s+(t-s))$ that corresponds to the predicted size and age at time $t$ of the dead individual.
\end{itemize}

\begin{definition}\label{def-eds}Let $Q(ds,di,d\theta)$ be a Poisson point measure (PPM) on $\R_+\times \N^*\times \R_+$ with intensity $ds\otimes n(di)\otimes d\theta$ where $ds$ and $d\theta$ are Lebesgue measures and where $n(di)$ is the counting measure on $\N^*$. The PPM provides possible times of events. For each time, we draw the label $i$ of the individual who may reproduce or die. The parameter $\theta$ allows to define whether a birth of death occurs. Assume also that the initial condition are $\bar{Z}_0$, $S_0$ and $V_0=\{1,\dots,I_0\}$, with $\E(\langle Z_0,1\rangle)<+\infty$ a.s. Then
\begin{multline}
\bar{Z}_t=  \sum_{i\in V_0}\delta_{(i,\Xi^i(t ; 0,\bar{Z}_0,S_0),a^i_0+t)}+\int_0^t \int_{\N^*\times \R_+} Q(ds,di,d\theta)\ \ind_{i\in V_{s_-}}\Big[ \label{sde}\\
\begin{aligned}
 & \Big(\delta_{(I_{s_-}+1,\Xi^{I_{s_-}+1}(t ; s,\bar{Z}_{s_-}+\delta_{(I_{s_-}+1,\xi_0,0)},S_s),t-s)}\\
 & \quad + \sum_{j\in V_{s_-}}\big(\delta_{(j,\Xi^j(t ; s,\bar{Z}_{s_-}+\delta_{(I_{s_-}+1,\xi_0,0)},S_s),a^j_s+(t-s))}-\delta_{(j,\Xi^j(t ; s,\bar{Z}_{s_-},S_s),a^j_s+(t-s))}\big)\Big)\\
 & \hspace{8cm}\times \ind_{\theta\leq m_1(i,s_-,\bar{Z}_{s_-},S_{s})}\\
 & \;+\Big(- \delta_{(i,\Xi^i(t ; s,\bar{Z}_{s_-},S_s),a^i_s+(t-s))}\\
& \quad +\sum_{j\in V_{s_-}}\big(\delta_{(j,\Xi^j(t ; s,\bar{Z}_{s_-}-\delta_{(i,\xi^i_{s_-},a^i_{s_-})},S_s),a^j_s+(t-s))} -\delta_{(j,\Xi^j(t ; s,\bar{Z}_{s_-},S_s),a^j_s+(t-s))}\big)\Big)\\
 & \hspace{6cm}\times \ind_{m_1(i,s_-,\bar{Z}_{s_-},S_{s})<\theta\leq m_2(i,s_-,\bar{Z}_{s_-},S_{s})}\Big]\\
\end{aligned}
\end{multline}where:
\begin{align}
m_1(i,s_-,\bar{Z}_{s_-},S_{s})= & \beta(\xi^i_{s_-},a^i_{s_-}, S_{s_-})\\
m_2(i,s_-,\bar{Z}_{s_-},S_{s})= & m_1(i,s_-,\bar{Z}_{s_-},S_{s_-})+\mu(\xi_{s},a^i_{s_-},S_{s_-}).
\end{align}
\end{definition}

\bigskip \medskip

\noindent In this definition the first of the two terms in the square brackets corresponds to the births and the second to the deaths. The first term starts with the addition of one new individual, followed by the corresponding updating of the future course of the trajectories of the other individuals in order to eventually get the right outcome at the final time $t$. The first part of the second term, dealing with the deaths,  starts with removing an individual, followed by a corresponding updating of the growth trajectories

Moment estimates obtained from \eqref{sde} are very important for the proofs. By adapting the proofs in \cite{champagnatferrieremeleard2} (Th. 2.5) (see also \cite{fourniermeleard}) where the main ingredient is the boundedness of the birth rate, we can show that:
\begin{lemma}Under Assumptions \ref{hypotheses-taux} and if $\E(\langle Z_0,1\rangle)<+\infty$ as in Def. \ref{def-eds}, then:\\
(i) For any $T>0$, $\E\big(\sup_{t\in [0,T]}\langle Z_t,1\rangle)<\E(\langle Z_0,1\rangle)e^{\bar{\beta}T}<+\infty$,\\
(ii) If we have additionally that $\E(\langle Z_0,1\rangle^p)<+\infty$, then
\begin{equation}
\E\big(\sup_{t\in [0,T]}\langle Z_t,1\rangle^p)<+\infty.\label{moment-p}
\end{equation}
\end{lemma}
As a consequence we obtain that (see \cite{champagnatferrieremeleard2,chithese}):
 \begin{proposition}
The process $(Z_t,S_t)_{t\in \R_+}$ is well and uniquely defined on $[0,T]$ for any initial condition such that $E(\langle Z_0,1\rangle)<+\infty$ and for any Poisson point measure $Q$. Moreover, the infinitesimal generator of $(Z_t)_{t\in \R_+}$ is \eqref{generateurinfinitesimal}.
 \end{proposition}

\subsection{Martingale problem}

We end this section with a martingale problem that will be useful to derive the large population limits. Heuristically, the decomposition of the process $\langle Z_t,f\rangle$, for any test function $f\in \Co^1$, into a predictable finite variation process and a square integrable martingale can be viewed as a description of the paths as solutions of the evolution equation associated with the generator $L$ (predictable finite variation part) plus noise (martingale part). The proof is given in Appendix \ref{AppendixA}.

\begin{proposition}\label{proposition:PBM}Let us assume that $\E(\langle Z_0,1\rangle^p)<+\infty$ for $p\geq 2$. Let us consider a test function $f(t,\xi,a)$ of class $\Co^1$. Then:
\begin{align}
M^f_t= & \langle Z_t,f(t,.,.)\rangle-\langle Z_0,f(0,.,.)\rangle-\int_0^t  \int_{\R_+^2} \Big(\frac{\partial f}{\partial s}(s,\xi,a)+\frac{\partial f}{\partial a}(s,\xi,a)+g(\xi,S_s)\frac{\partial f}{\partial \xi}(s,\xi,a)\nonumber\\
& +  f(s,\xi_0,0)\beta(\xi,a ,S_s)-f(s,\xi,a)\mu(\xi,a,S_s)\Big)Z_s(d\xi,da)ds,\label{pbm}
\end{align}is a square integrable martingale with predictable quadratic variation:
\begin{align}
\langle M^f\rangle_t= & \int_0^t \int_{\R_+^2}\big( f^2(s,\xi_0,0) \beta(\xi,a,S_s)+f^2(s,\xi,a)\mu(\xi,a,S_s)\big) Z_s(d\xi,da)\,ds.\label{bracket}
\end{align}
\end{proposition}

\section{Large populations}

We now focus on large populations of Daphnia. To that end we consider a sequence of processes describing the evolution of the population through time when starting with an initial condition of size proportional to the integer parameter $K$ \eqref{food} that we let increase to infinity. As already indicated, we may think of this $K$ as the volume in which the population and its food live. We moreover scale the population with ${1}\over{K}$, i.e., we transform from population size to population density.
We thus consider a sequence $\bar{Z}^K$ of populations such that:
\begin{equation}
\bar{Z}^K_t(dj,d\xi,da)=\frac{1}{K}\sum_{i\in V_t^K} \delta_{(i,\xi^i_t,a^i_t)}(dj,d\xi,da)
\end{equation}where $V^K_t$ is the set of individuals alive at time $t$. Again, $Z^K_t(d\xi,da)=\bar{Z}^K_t(\N^*\times d\xi\times da)$ is the marginal on $\R_+^2$ of $\bar{Z}^K_t$. We also consider $(S^K_t)_{t\in \R_+}$ the sequence of associated food concentrations, also indexed by $K$ and satisfying:
\begin{equation}
\frac{dS^K_t}{dt}=f(S^K_t)-\int_{\N^*\times \R_+^2} \gamma(\xi,S^K_t)\bar{Z}^K_t(dj,d\xi,da),
\end{equation}
with initial conditions $(S^K_0)_{K\in \N^*}$ that converge in probability to $S_0\in \R_+$. \\
For each $K\in \N^*$, $K\bar{Z}^K$ has the dynamics of the process introduced in Definition \ref{def-eds} with initial conditions $K\bar{Z}^K_0$ for which we assume that
\begin{equation}
\sup_{K\in \N^*} \E\Big(\langle Z^K_0, 1\rangle^2 \Big)<+\infty.\label{hypotheses-moments}
\end{equation}

\begin{proposition}\label{prop-convergence}Let $T>0$. Under the Assumptions \ref{hypotheses-taux}, the sequence $(Z^K,S^K)_{n\geq 1}$ introduced in this section converges in probability in $\D([0,T],\mathcal{M}_F(\R_+^2)\times \R_+)$ to the unique continuous solution $(\zeta,\varrho)$ of the following deterministic equations, characterized for any function $f(t,\xi,a)$ of class $\Co^1$ by:
\begin{align}
& \langle \zeta_t,f(t,.,.)\rangle=  \langle \zeta_0,f(0,.,.)\rangle+ \int_0^t \Big\langle \zeta_s,\frac{\partial f(s,.,.)}{\partial s}+\frac{\partial f(s,.,.)}{\partial a}+g \frac{\partial f(s,.,.)}{\partial \xi}\Big\rangle ds\nonumber\\
& \hspace{2cm} +  \int_0^t \langle \zeta_s,f(s,\xi_0,0)\beta(.,.,\varrho_s)-f(s,.,.)\mu(.,.,\varrho_s)\rangle ds.\label{pde}\\
& \frac{d\varrho}{dt}= f(\varrho(t))-\int_{\R_+^2} \gamma(\xi,\varrho(t))\zeta_t(d\xi,da).\label{pde2}
\end{align}
\end{proposition}
The proof of Proposition \ref{prop-convergence} is given in Appendix \ref{AppendixB}.\\

We conclude by showing that equations \eqref{pde}-\eqref{pde2} allow us to recover the equations of \cite{diekmanngyllenbergmetznakaokaderoos}. We first establish a precise form of the solution $\zeta_t$. That this measure is not absolutely continuous with respect to the Lebesgue measure on $\R_+^2$, even if the initial condition $\zeta_0$ is, was already noticed in e.g. \cite{chithese,trangdesdev}. There it was stated that age and size where both parameterized by time, but no precise form for the measure for these cases was given.

\begin{corollary}\label{corollaire}Assume that the initial condition has a marginal measure in age that is absolutely continuous with respect to the Lebesgue measure on $\R_+$ so that $\zeta_0(d\xi,da)= q_0(a,d\xi) da$, where $q_0(a,d\xi)$ is a transition measure on $\R_+$. We denote by $\nu_{a,t}(d\xi)$ the image measure of $q_0(a-t,d\xi)$ through the application $\xi\mapsto \Xi(t ; 0, \xi,\zeta_{0},\varrho_{0})$.
For any $t\in \R_+$, the marginal $\zeta_t$ at time $t$ of the solution of \eqref{pde} is a.s. given by:
\begin{align}
\zeta_t(d\xi,da)= & \ind_{a< t} \ b(t-a)\mathcal{F}(a,t-a,\varrho_{[0,t]}) \delta_{\Xi(t ; t-a,0,\zeta_{t-a},\varrho_{t-a})}(d\xi)\ da \nonumber\\
+ & \ind_{a\geq t}\  \mathcal{F}'(t,a-t,\xi', \varrho_{[0,t]})\ind_{\xi=\Xi(t; 0,\xi',\zeta_0,\varrho_0)} \nu_{a,t}(d\xi)\ da\label{formezeta}
\end{align}where $b(t)=\int_{\R_+^2} \beta(\xi,a,\varrho_{t})\zeta_t(d\xi,da)$ is the total birth rate at time $t$, where $\varrho_{[0,t]}=(\varrho_s)_{s\in [0,t]}$ and where for $a< t$,
\begin{equation}\mathcal{F}(a,t_0,\varrho_{[0,t_0+a]})=\exp\Big(-\int_0^a \mu\big(\Xi(t_0+\alpha ; t_0,\xi_0,\zeta_{t_0},\varrho_{t_0}),\alpha ,\varrho_{t_0+\alpha}\big)d\alpha\Big)\label{survie-corol}\end{equation}is the probability that an individual born at $t_0$ survives until age $a$ when the food environment is given by $(\varrho_s)_{s\in [0,t_0+a]}$. For $a\geq t$,
\begin{equation}\mathcal{F}'(t,a_0,\xi',\varrho_{[0,t]})=\exp\Big(-\int_{0}^t \mu\big(\Xi(s ; 0,\xi',\zeta_{0},\varrho_{0}),a_0+s ,\varrho_{s}\big)ds\Big)\label{survie-corol2}\end{equation}is the probability that an individual alive at $t=0$ with age $a_0$ and size $\xi'$ survives until time $t$ in an environment $\varrho_{[0,t]}$.
\end{corollary}

With Corollary \ref{corollaire}, we recover the equations of  \cite{diekmanngyllenbergmetznakaokaderoos}. Equation \eqref{growth} provides the deterministic differential equation describing the growth of Daphnias, represented by the distribution $\zeta_t(d\xi,da)$ (see \eqref{formezeta}):
\begin{align*}
 \frac{d\xi}{da}(a)=g\big(\xi(a),\varrho(t_0+a)\big); \quad \xi(0)=\xi_0.
\end{align*}
If we consider an individual born at time $t_0>0$ and follow the survival probability through time
$a\mapsto \mathcal{F}(a,t_0,\varrho_{[0,t_0+a]})$, Equation \eqref{survie-corol} gives the decay of the survival probability of an individual of age $a$ at time $t$:
\begin{equation*}
\frac{d\mathcal{F}}{da}(a,t_0,\varrho_{[0,t_0+a]})=-\mu\big(\Xi(t_0+a ; t_0,\xi_0,\zeta_{t_0},\varrho_{t_0}),a,\varrho_{t_0+a}\big)\mathcal{F}(a,t_0,\varrho_{[0,t_0+a]}).
\end{equation*}
From \eqref{formezeta}, we have the Daphnia population birth rate at time $t$:
\begin{align*}
& b(t)=\int_0^t \beta(\Xi(t;t-a,0,\zeta_{t-a},\varrho_{t-a}),a,\varrho_t)b(t-a)\mathcal{F}(a,t-a,\varrho_{[0,t]}) da \\
& \hspace{1cm} + \int_{\R_+^2} \beta(\Xi(t;0,\xi,\zeta_0,\varrho_0),a+t,\varrho_t)\mathcal{F}'(t,a,\xi,\varrho_{[0,t]}) \zeta_0(d\xi,da)
\end{align*}
The first term represents the contributions of individuals born after time $0$, while the second term corresponds to individuals who where present initially. From this, we can deduce the algal concentration:
\begin{align*}
& \frac{d\varrho}{dt}(t)=f(\varrho(t))-\int_0^t \gamma(\Xi(t;t-a,0,\zeta_{t-a},\varrho_{t-a}),\varrho(t)) b(t-a)\mathcal{F}(a,t-a,\varrho_{[0,t]}) da\\
& \hspace{1cm} + \int_{\R_+^2} \gamma(\Xi(t;0,\xi,\zeta_0,\varrho_0),\varrho(t))\mathcal{F}'(t,a,\xi,\varrho_{[0,t]})\zeta_0(d\xi,da).
\end{align*}

\begin{proof} {\it of Corollary \ref{corollaire}}
 For the proof, we start by showing that $\zeta_t(d\xi,da)$ admits a density $m(\xi,a,t)$ w.r.t. a dominating measure underlying \eqref{formezeta}. The equations satisfied by $m(\xi,a,t)$ are then derived by separating the domain into $\R_+\times \{a\geq t\}$ and $\R_+\times \{a<t\}$, which corresponds to first studying the individuals born before and after initial time.

First, recall that there is a unique solution $\varrho$ to \eqref{pde2}. \\
Let $\varphi\in \Co^1_{\rm{b}}(\R_+^2,\R)$, let $t\in \R_+$ and  consider the associated test function:
$$f(s,\xi,a)=\varphi\big(\Xi(t ; s, \xi,\zeta_{s},\varrho_{s}),a+t-s\big).$$
This function $f$ is the unique solution of:
\begin{align*}
\Big( \frac{\partial f}{\partial s}+\frac{\partial f}{\partial a}+g \frac{\partial f}{\partial \xi} \Big)(s,\xi,a)=0,\qquad f(t,\xi,a)=\varphi(\xi,a)
\end{align*}(e.g. \cite{evans}). As a consequence, using this test function $f$ in \eqref{pde}:
\begin{multline}
\langle \zeta_t,\varphi\rangle
=  \langle \zeta_0,\varphi(\Xi(t ; 0, \xi,\zeta_{0},\varrho_{0}),.+t)\rangle\\
+\int_0^t \Big(\varphi(\Xi(t ; s, \xi_0,\zeta_{s},\varrho_{s}),t-s)\langle \zeta_s,\beta(.,.,\varrho_s)\rangle \\
- \int_{\R_+^2} \varphi\big(\Xi(t ; s, \xi,\zeta_{s},\varrho_{s}),a+t-s\big)\mu(\xi,a,\varrho_s)\zeta_s(d\xi,da)\Big) ds.\label{etape2}\end{multline}The first term is related to individuals that are alive at time 0. The second integral relates to births between time 0 and time $t$. The third term corresponds to the deaths between time 0 and time $t$. \\
If we consider positive functions $\varphi$, then, neglecting the non-positive terms in \eqref{etape2}:
\begin{align}0\leq \langle \zeta_t,\varphi\rangle \leq & \int_{t}^{+\infty} \Big(\int_{\R_+} \varphi(\Xi(t ; 0, \xi,\zeta_{0},\varrho_{0}),a) q_0(a-t,d\xi)\Big) da \nonumber\\
 & \hspace{1cm}+\int_0^t b(t-a)\varphi(\Xi(t ; t-a,\xi_0,\zeta_{t-a},\varrho_{t-a}),a) da.
\end{align}Notice that the population is naturally divided into two sets. Since the aging velocity is 1, the individuals who were alive at initial time are of age greater than $t$ at time $t$. Individuals born after time 0 are of age smaller than $t$. So, if the function $\varphi$ has support included in the set $\R_+\times \{a< t\}$, then we see that on $\R_+\times \{a< t\}$, $\zeta_t(d\xi,da)$ is absolutely continuous with respect to $\delta_{\Xi(t ; t-a,\xi_0,\zeta_{t-a},\varrho_{t-a})}(d\xi)\, da$. Similarly, on the set $\R_+\times \{a\geq t\}$, $\zeta_t(d\xi,da)$ admits a density with respect to $\nu_{a,t}(d\xi)\, da$. Denote by $m(\xi,a,t)$ the density of $\zeta_t$ with respect to the measure $\ind_{a<t}\, \delta_{\Xi(t ; t-a,\xi_0,\zeta_{t-a},\varrho_{t-a})}(d\xi)\, da+\ind_{a\geq t} \nu_{a,t}(d\xi)\, da$. \\
Substituting this density in the third term of \eqref{etape2} gives, for $\varphi$ with support in $\R_+\times \{a\geq t\}$:
\begin{align*}
\langle \zeta_t,\varphi\rangle= & \int_{t}^{+\infty} \int_{\R_+} \varphi(\Xi(t ; 0, \xi,\zeta_{0},\varrho_{0}),a) q_0(a-t,d\xi) da \nonumber\\
 -& \int_t^{+\infty} da \int_{\R_+}q_0(a-t,d\xi)\ \varphi\big(\Xi(t;0,\xi,\zeta_0,\varrho_0),a\big)\nonumber\\
  & \hspace{1cm}\times \int_0^{t} \Big[\mu(\Xi(s;0,\xi,\zeta_0,\varrho_0),a-t+s),\varrho_s)\nonumber\\
  & \hspace{2.5cm} m(\Xi(s;0,\xi,\zeta_0,\varrho_0),a-t+s,s) \Big]ds
\end{align*}By identification, the density $m(\xi,a,t)$ of $\zeta_t$, on $\{a\geq t\}$, w.r.t. $\nu_{a,t}(d\xi) da$ satisfies:
\begin{align*}
m(\xi,a,t)= & 1-\int_0^t \mu(\Xi(s;0,\xi',\zeta_0,\varrho_0),a-t+s),\varrho_s)\nonumber\\
& \hspace{1.5cm} m(\Xi(s;0,\xi',\zeta_0,\varrho_0),a-t+s,s)\ ds\ \ind_{\{\xi=\Xi(t;0,\xi',\zeta_0,\varrho_0)\}}
\end{align*}where there exists a unique $\xi'$ such that $\xi=\Xi(t;0,\xi',\zeta_0,\varrho_0)$ under Assumptions \ref{hypotheses-taux}. Notice that $a-t$ is the age of the individual at time $0$ and is not a real function of time. Thus, we recognize an ordinary differential equation of degree 1 for $s\mapsto m(\Xi(s;0,\xi',\zeta_0,\varrho_0),a-t+s,s)$ from which
\begin{align*}
m(\xi,a,t) = & \exp\Big(-\int_0^t \mu(\Xi(s;0,\xi',\zeta_0,\varrho_0),a-t+s,\varrho_s) ds\Big)\ind_{\{\xi=\Xi(t;0,\xi',\zeta_0,\varrho_0)\}}.
\end{align*}
This yields the second part of \eqref{formezeta}.

Choosing $\varphi$ with support in $\R_+\times \{a<t\}$, \eqref{etape2} yields:
\begin{align*}
 \langle \zeta_t,\varphi\rangle =&\int_0^t da\ b(t-a)\varphi\big(\Xi(t ; t-a,\xi_0,\zeta_{t-a},\varrho_{t-a}),a\big) \nonumber\\
 -  & \int_0^{t} da \ \varphi(\Xi(t;t-a,\xi_0,\zeta_{t-a},\rho_{t-a}),a) \nonumber\\
 &\hspace{1cm}\times  \int_{0}^a \Big[\mu\big(\Xi(u+t-a;t-a,\xi_0,\zeta_{t-a},\rho_{t-a}),u,\rho_{u+t-a}\big)\nonumber\\
  & \hspace{2.5cm} m\big(\Xi(u+t-a;t-a,\xi_0,\zeta_{t-a},\rho_{t-a}),u,u+t-a\big) \Big]du\nonumber\\
  = & \int_0^t da\ b(t-a)\varphi\big(\Xi(t ; t-a,\xi_0,\zeta_{t-a},\varrho_{t-a}),a\big) \nonumber\\
 -  & \int_0^{t} da \ \varphi(\Xi(t;t-a,\xi_0,\zeta_{t-a},\rho_{t-a}),a) \nonumber\\
 &\hspace{1cm}\times  \int_{0}^a\int_{\R_+} \mu\big(\xi,u,\rho_{u+t-a}\big) \zeta_{u+t-a}(d\xi,du)
\end{align*}Thus, on $\{a<t\}$, the density $m(\xi,a,t)$ of $\zeta_t$ w.r.t. $\delta_{\Xi(t ; t-a,\xi_0,\zeta_{t-a},\varrho_{t-a})}(d\xi)\, da$ satisfies
\begin{align*}
m(\xi,a,t)=& b(t-a) -\int_{0}^a\int_{\R_+} \mu\big(\xi',u,\rho_{u+t-a}\big) \zeta_{u+t-a}(d\xi',du).
\end{align*}Notice that $t-a$ is the time of birth of the individual and is not a real function of age. Then we recognize again an ordinary differential equation of order 1 from which
\begin{align*}
m(\xi,a,t)= & b(t-a) \exp\Big(-\int_0^a \mu\big(\Xi(u+t-a;t-a,\xi_0,\zeta_{t-a},\rho_{t-a}),u,\rho_{u+t-a}\big) du \Big).
\end{align*}This ends the proof of the announced result \eqref{formezeta} for $a<t$.
\end{proof}

\section{Conclusion}

In this paper we rigorously underpin the long suspected individual-based nature of physiological structured population models as studied by the team of researchers surrounding Odo Diekmann. Such results were already known for purely age-based models (e.g. \cite{oelschlager,chithese,trangdesdev,ferrieretran,jagersklebaner,jagersklebaner2}), and our present result still applies only to a subclass of the models considered by Odo and his co-workers,  although a paradigmatic one, to wit their so-called ''Daphnia'' models e.g. \cite{diekmanngyllenbergmetznakaokaderoos,DiekmannMetz2010}. The i-state variables of these models move in a continuous deterministic fashion, dependent on their own value and the environment (food: ''algae''), and births occur to a single fixed birth state at a rate that depends on the i-state and the environment (and are not e.g. coupled to specific i-state transitions). Even for this  restricted class of physiologically structured population models some twiddling of the existing probabilistic toolbox was in order. Yet, overall the toolbox proved its mettle, and the long standing assumption of a solid individual-based foundation of the theory of physiologically structured populations was duly vindicated. The longer term goal is to extend this vindication to the full class of models put forward in \cite{MetzDiekmann,Diekmannetal1998,Diekmannetal2001,Diekmannetal2003}.

\bigskip

\noindent
\textbf{Acknowledgements}: This work benefitted from the support from the ``Chaire Mod\'elisation Math\'ematique et Biodiversit\'e of Veolia Environnement - Ecole Polytechnique - Museum National d'Histoire Naturelle - Fondation X".

\bigskip

\appendix

\begin{flushleft} {\large{ \textbf {Appendices}}} \end{flushleft}

\section{Proof of Proposition \ref{proposition:PBM} }\label{AppendixA}
Let $f(t,\xi,a)$ be a function of class $\Co^1$. From \eqref{sde}, we obtain
\begin{align}
\lefteqn{\langle Z_t,f\rangle=  \sum_{i\in V_0}f(t,\Xi^i(t ; 0,\bar{Z}_0,S_0),a^i_0+t)+\int_0^t \int_{\N^*\times \R_+} Q(ds,di,d\theta)\ \ind_{i\in V_{s_-}}\Big[}\nonumber\\
 & \ind_{\theta\leq m_1(i,s_-,\bar{Z}_{s_-},S_{s})} \Big(f(t,\Xi(t ; s,0,\bar{Z}_{s}+\delta_{(I_{s_-}+1,\xi_0,0)},S_s),t-s)\nonumber\\
 & \quad + \sum_{j\in V_{s_-}}\big(f(t,\Xi^j(t ; s,\bar{Z}_{s_-}+\delta_{(I_{s_-}+1,\xi_0,0)},S_s),a^j_s+(t-s))\nonumber\\
 & \hspace{5cm}-f(t,\Xi^j(t ; s,\bar{Z}_{s_-},S_s),a^j_s+(t-s))\big)\Big)\nonumber\\
+ & \ind_{m_1(i,s_-,\bar{Z}_{s_-},S_s)<\theta\leq m_2(i,s_-,\bar{Z}_{s_-},S_{s})}\Big(- f(t,\Xi^i(t ; s,\bar{Z}_{s_-},S_s),a^i_s+(t-s))\nonumber\\
 & \quad +\sum_{j\in V_{s_-}}\big(f(t,\Xi^j(t ; s,\bar{Z}_{s_-}-\delta_{(i,\xi^i_{s_-},a^i_{s_-})},S_s),a^j_s+(t-s))\nonumber\\
 & \hspace{5cm}-f(t,\Xi^j(t ; s, \bar{Z}_{s_-},S_s),a^j_s+(t-s))\big) \Big)\Big].\label{l0}
\end{align}
Using \eqref{growth}, we have for any $s<t$:
\begin{align*}
f(t,\Xi^i(t ; s,\bar{Z}_s,S_s),a^i_s+(t-s))
= & f(s,\xi^i_s,a^i_s)+\int_s^t \Big(\frac{\partial f}{\partial u}+\frac{\partial f}{\partial a}(u,\Xi^i(u ; s,\bar{Z}_s,S_s),a^i_s+u-s)\nonumber\\
+ & g(\Xi^i(u ; s,\bar{Z}_s,S_s),S_u)\frac{\partial f}{\partial x}(u,\Xi^i(u ; s,\bar{Z}_s,S_s),a^i_s+u-s)\Big)du
\end{align*}
Recall that we denoted by $T_k$, $k\geq 1$ the birth and death events in the population. By convention, we let $T_0=0$. Let us consider an individual $i$. Let $t_0\in \{T_k, k\geq 0\}$ be the birth time of the individual (or 0 if the individual is alive at time $0$) and $a^i_{t_0}$ be its age at time $t_0$ (0 if $t_0$ is the birth time). The sum of the terms in the r.h.s. of \eqref{l0} associated with individual $i$ is equal to:
\begin{multline*}
 f(t_0,\xi^i_{t_0},a^i_{t_0})+\sum_{k\geq 0} \int_{t\wedge T_k \vee t_0}^{t\wedge T_{k+1} \vee t_0} \Big(\frac{\partial f}{\partial u}+\frac{\partial f}{\partial a}(s,\Xi^i(s ; T_k,\bar{Z}_{T_k},S_{T_k}),a^i_{t_0}+s-t_0)\\
+  g(s,\Xi^i(s ; T_k,\bar{Z}_{T_k},S_{T_k}),S_s)\frac{\partial f}{\partial x}(s,\Xi^i(s ; T_k,\bar{Z}_{T_k},S_{T_k}),a^i_{t_0}+s-t_0)\Big)ds\\
- \int_0^t \int_{\N^*\times \R_+}  \ind_{j=i ; i\in V_{s_-}}\ind_{m_1(i,s_-,\bar{Z}_{s_-},S_{s})<\theta\leq m_2(i,s_-,\bar{Z}_{s_-},S_{s})}f(s,\xi^i_{s_-},a^i_{s_-})\Big]dQ.
\end{multline*}The last integral correspond to the death term when individual $i$ is dead before $t$. Thus, \eqref{l0} gives:
\begin{align*}
\lefteqn{\langle Z_t,f(t,.,.)\rangle =  \sum_{i\in V_0} \Big[f(0,\xi^i_{0},a^i_{0})}\\
+ & \sum_{k\geq 0} \int_{t\wedge T_k}^{t\wedge T_{k+1}}  \Big(\frac{\partial f}{\partial u}+\frac{\partial f}{\partial a}(s,\Xi^i(s ; T_k,\bar{Z}_{T_k},S_{T_k}),a^i_{0}+s)\\
 &\hspace{2cm}+ g(s,\Xi^i(s ; T_k,\bar{Z}_{T_k},S_{T_k}),S_s)\frac{\partial f}{\partial x}(s,\Xi^i(s ; T_k,\bar{Z}_{T_k},S_{T_k}),a^i_{0}+s)\Big)ds\\
- & \int_0^t \int_{\N^*\times \R_+} Q(ds,dj,d\theta) \ind_{j=i}\ind_{i\in V_{s_-}}\ind_{m_1(i,s_-,\bar{Z}_{s_-},S_{s})<\theta\leq m_2(i,s_-,\bar{Z}_{s_-},S_{s})}f(s,\xi^i_{s_-},a^i_{s_-})\Big]\\
+ & \int_0^t \int_{\N^*\times \R_+} Q(ds,di,d\theta)\ \ind_{i\in V_{s_-}\setminus V_0}\Big[\\
 & \Big(f(s,\xi_0,0)+\sum_{k\geq 0} \int_{t\wedge T_k\vee s}^{t\wedge T_{k+1}\vee s}  \Big(\frac{\partial f}{\partial u}+\frac{\partial f}{\partial a}(u,\Xi^{I_{s_-}+1}(u ; T_k,\bar{Z}_{T_k},S_{T_k}),u-s)\\
+ & g(u,\Xi^{I_{s_-}+1}(u ; T_k,\bar{Z}_{T_k},S_{T_k}),S_u)\frac{\partial f}{\partial x}(u,\Xi^{I_{s_-}+1}(u ; T_k,\bar{Z}_{T_k},S_{T_k}),u-s)\Big)du\Big)\\
& \hspace{2cm}\ind_{\theta\leq m_1(i,s_-,\bar{Z}_{s_-},S_{s})}\\
- & f(s,\xi^i_{s_-},a)\ind_{m_1(i,s_-,\bar{Z}_{s_-},S_{s})<\theta\leq m_2(i,s_-,\bar{Z}_{s_-},S_{s})}\Big],
\end{align*}where the first bracket corresponds to individuals alive at time $0$ and where the second bracket correspond to individuals born after time $0$.
For $s<u$
$$\sum_{i\in V_{s_-}} \delta_{(\Xi^i(u ; s,\bar{Z}_s,S_s),a^i_s+u-s)}(d\xi,da)=Z_u(d\xi,da)$$
provided there has been no jumps between $s$ and $u$. Thus, we have:
\begin{align*}
\langle Z_t,f(t,.,.)\rangle   
=& \langle Z_0,f(0,.,.)\rangle+ \int_0^t ds \int_{\R_+^2}Z_s(d\xi,da)  \Big(\frac{\partial f}{\partial s}+\frac{\partial f}{\partial a}(s,\xi,a)+ g(\xi,S_s)\frac{\partial f}{\partial \xi}(s,\xi,a)\Big)\\
+& \int_0^t \int_{\R_+} \big( f(s,\xi_0,0) \beta(\xi,a,S_s)-f(s,\xi,a)\mu(\xi,a,S_s)\big)Z_s(d\xi,da)\\
+ & \int_0^t \int_{\N^*\times \R_+}\ind_{i\in V_{s_-}}\Big[f(s,\xi_0,0)\ind_{\theta\leq m_1(i,s_-,\bar{Z}_{s_-},S_{s})}\\
-  &  f(s,\xi^i_{s_-},a)\ind_{m_1(i,s_-,\bar{Z}_{s_-},S_{s})<\theta\leq m_2(i,s_-,\bar{Z}_{s_-},S_{s})}\Big] \widetilde{Q}(ds,di,d\theta),
\end{align*}where $\widetilde{Q}(ds,di,d\theta)=Q(ds,di,d\theta)-ds\otimes n(di)\otimes d\theta$ is the compensated Poisson point measure associated with $Q$. The integral with respect to (w.r.t.) $\widetilde{Q}(ds,di,d\theta)$ provides the martingale $M^f$. This achieves the proof.\hfill $\Box$

\section {Sketch of the proof of Proposition \ref{prop-convergence}}\label{AppendixB}
When starting from \eqref{pbm} and using controls of moments as in \cite{fourniermeleard}, the proof is similar to the one in \cite{trangdesdev,chithese}.\\

\noindent \textbf{Step 1} We start by noticing that under the Assumption \eqref{hypotheses-moments}, we have the following estimate (e.g. \cite{champagnatferrieremeleard2}):
\begin{equation}
\sup_{n\in \N^*} \E\Big(\sup_{t\in [0,T]} \langle Z^n_t, 1\rangle^2 \Big)<+\infty.\label{controlemoments}
\end{equation}Moreover, from Assumptions \ref{hypotheses-taux} and \eqref{eq1}, the size of any individual is bounded on $[0,T]$ by $\bar{\xi}=\xi_0+\bar{g}T$ and there exists for every $\varepsilon$ a non random constant $\bar{S}_\varepsilon$ such that:
\begin{equation}
\sup_{n\in \N^*}\P\Big(\sup_{t\in [0,T]}S^n_t>\bar{S}_\varepsilon\Big)<\varepsilon.
\end{equation}From these estimates and Assumption \ref{hypotheses-taux} (ii), there exists a constant $A_\varepsilon \in (0,A)$ such that:
\begin{equation}
\sup_{n\in \N^*}\P\Big(\sup_{t\in [0,T]}Z^n_t\big([\xi_0,\bar{\xi}]\times [0,A_\varepsilon]\big)>\varepsilon\Big)<\varepsilon.
\end{equation}

\noindent \textbf{Step 2} It is easy to see that the limiting values of $(Z^n,S^n)_{n\in \N^*}$ are necessary continuous. Let us check the $C$-tightness (e.g. \cite{jacod}) of $(Z^n,S^n)_{n\in \N^*}$ in $\D([0,T],\mathcal{M}_F(\R_+^2)\times \R_+)$. Using a criterion by \cite{meleardroelly} and given the compact containment that follows from Step 1, it is sufficient to prove the tightness of $(S^n)_{n\in \N^*}$ and of the predictable finite variation part and martingale part of $(\langle Z^n,f\rangle)_{n\in \N^*}$ for $f$ in $\Co_b^1(\R_+^3,\R)$ (which contains the constant function equal to 1). This is obtained by using Aldous-Rebolledo criteria (e.g. \cite{joffemetivier}) and adapting for instance \cite{champagnatferrieremeleard2,trangdesdev} with the estimates of Step 1. \\

\noindent \textbf{Step 3} The identification of the martingale problem satisfied by the limiting values provides \eqref{pde}-\eqref{pde2}. Uniqueness of the solution of \eqref{pde}-\eqref{pde2} stems from the Assumptions \ref{hypotheses-taux}. As a consequence, there is a unique limiting value and we have convergence in distribution of $(Z^n,S^n)_{n\in \N^*}$ to the solution $(\zeta,\varrho)$. Since the latter is deterministic, the convergence is also a convergence in distribution.



{\footnotesize
\providecommand{\noopsort}[1]{}\providecommand{\noopsort}[1]{}\providecommand{\noopsort}[1]{}\providecommand{\noopsort}[1]{}

}
\end{document}